\documentclass[11pt,twoside,leqno]{aomamlt2e}
  \pageno{1}
\received{}
\revised{}

\newtheorem{theorem}{Theorem}[section]
\newtheorem{lemma}{Lemma}[section]

\begin{document}
\currannalsline{0}{2006}

\title{A New Upper Bound for Diagonal Ramsey Numbers}

\acknowledgements{The author is kindly supported by a grant from St
John's College, Cambridge.}
\author{David Conlon}

\institution{Department of Pure Mathematics and Mathematical
Statistics, University of Cambridge, Cambridge, United Kingdom\\
\email{D.Conlon@dpmms.cam.ac.uk}}

\shorttitle{Diagonal Ramsey Numbers}

\begin{abstract}
We prove a new upper bound for diagonal two-colour Ramsey numbers, 
showing that there exists a constant $C$ such that
\[r(k+1, k+1) \leq k^{- C \frac{\log k}{\log \log k}} \binom{2k}{k}.\]
\end{abstract}

\section{Introduction}

The Ramsey number $r(k,l)$ is the smallest natural number $n$ such that, 
in any red and blue colouring of the edges of the complete graph on $n$ 
vertices, we are guaranteed to find either a red $K_k$ or a blue $K_l$.

That these numbers exist is a consequence of Ramsey's original theorem 
\cite{R30}, but the standard upper bound
\[r(k+1, l+1) \leq \binom{k+l}{k}\]
is due to Erd\H{o}s and Szekeres \cite{ES35}. 

Very little progress was made on improving this bound until the mid-eighties,
when a number of successive improvements were given, 
showing that, as expected, $r(k+1, l+1) = o(\binom{k+l}{k})$. 
Firstly, R\"{o}dl showed that for some constants $c$ and $c'$ we have
\[r(k+1, l+1) \leq \frac{c \binom{k+l}{k}}{\log^{c'}(k+l)}.\]
This result was never published, but a weaker bound,
\[r(k+1, l+1) \leq \frac{6 \binom{k+l}{k}}{\log \log(k+l)}\]
appears in the survey paper concerning Ramsey bounds by Graham and 
R\"{o}dl \cite{GR87}.    

Not long after these bounds were proven, Thomason \cite{T88} 
proved that there was a positive constant $A$ such that, for $k \geq l$, 
\[r(k+1, l+1) \leq \exp \left\{ -\frac{l}{2k} \mbox{log} k + A 
\sqrt{\log k} \right\} \binom{k+l}{k},\]
this being a major improvement on R\"{o}dl's bound when $k$ and $l$ are of 
approximately the same order, implying in particular that
\[r(k+1, k+1) \leq k^{-1/2 + A/\sqrt{\log k}} \binom{2k}{k}.\]

In this paper we will show how to improve on Thomason's result, obtaining

\begin{theorem}
There exists a constant $C$ such that 
\[r(k+1, k+1) \leq k^{- C \frac{\log k}{\log \log k}} \binom{2k}{k}.\]
\end{theorem}

In particular we have the following natural extension of Thomason's 
Theorem:

\begin{theorem} 
For all $s>0$ there exists a constant $C_s$ such that
\[r(k+1, k+1) \leq \frac{C_s}{k^s} \binom{2k}{k}.\]
\end{theorem}

\section{An Outline of the Proof}

Our argument (and also Thomason's) begins by assuming that we are trying 
to prove a bound of the form $r(k+1, l+1) \leq f(k, l) \binom{k+l}{k}$, 
where $f(k,l)$ is some slowly changing function in $k$ and $l$. In order 
to construct an inductive argument we will assume that for some such 
function we have $r(a+1, b+1) \leq f(a,b) \binom{a+b}{a}$ whenever $a$ is 
less than $k$ or $b$ is less than $l$, and that we would like to show 
that the same holds for $a = k$ and $b = l$.

To this end, let us suppose that $n = \lfloor f(k, l) \binom{k+l}{k} 
\rfloor = f^*(k,l) \binom{k+l}{k}$, say. Then by the argument that proves 
the Erd\H{o}s-Szekeres inequality
\[ r(k+1, l+1) \leq r(k, l+1) + r(k+1,l),\]
we see that, within any red/blue colouring of the edges of $K_n$ which 
does not contain a red $K_{k+1}$ or a blue $K_{l+1}$, every vertex $x$ 
can have red degree at most $r(k, l+1) - 1$ and blue degree at most 
$r(k+1, l) - 1$. Therefore, if $d_x$ is the red degree of the vertex $x$ 
(so that $n-1-d_x$ is the blue degree),
\begin{eqnarray*}
d_x & < & r(k, l+1)\\
& \leq & f(k-1, l) \binom{k+l-1}{k-1} \\
& = & \frac{f(k-1, l)}{f^*(k,l)} \frac{k}{k+l} n.
\end{eqnarray*}
Similarly, we may use the fact that $n - 1 - d_x \leq r(k+1, l) - 1$ to 
show that 
\[d_x \geq \left( 1 - \frac{f(k, l-1)}{f^*(k,l)} \frac{l}{k+l} \right) n.\]
Now, note that if $f$ were always one, then we would know that $d_x$ was 
less than $\frac{k}{k+l} n$ for each vertex $x$ and also that it was 
greater than or equal to $\frac{k}{k+l} n$, a contradiction which is 
equivalent to the Erd\H{o}s-Szekeres argument.

If, instead, we allow the size of $f(k, l)$ to change with both $k$ and 
$l$, albeit slowly, then we find that for each vertex $x$ the red degree 
$d_x$ is not much greater than $\frac{k}{k+l}n$ nor much less than it. So 
we find that the graph is approximately regular in degree, the proximity 
to true regularity being dependent upon how slowly $f(k,l)$ changes.

This approximate degree-regularity is not however the only structural 
information that we have about graph colourings which contain neither a 
red $K_{k+1}$ nor a blue $K_{l+1}$. We also know, for example, that in 
such a graph any red edge can lie in at most $r(k-1, l+1) - 1$ red 
triangles, and if the vertices of this red edge are $x$ and $y$ 
then there are at most $r(k,l) - 1$ vertices which are connected to $x$ 
by a red edge and $y$ by a blue edge. If we let $d_{x y}$ be the number 
of vertices which are connected to both $x$ and $y$ by a red edge, then 
these two conditions are enough to tell us that
\[d_{x y} \approx \left( \frac{k}{k+l} \right)^2 n,\]
the exact proximity being again dependent upon the rate at which $f(k,l)$ 
changes. That is, providing that we don't try and improve too much on the 
Erd\H{o}s-Szekeres bound, we can conclude that across any red edge we 
have approximately the expected number of red triangles (that would be in 
a random graph formed by choosing red edges with probability 
$\frac{k}{k+l}$). As a consequence, we see that across any red edge there 
are approximately the expected number of red $C_4$s of which the red edge 
is a diagonal. Importantly, this latter result is not restricted to red 
edges alone - it is straightforward to use the degree-regularity 
conditions and the analogous condition that we have approximately the 
expected number of blue $C_4$s across a blue edge in order to show that 
we have approximately the expected number of red $C_4$s across that edge 
as well. 

At this stage it is appropriate to recall (at least roughly) the 
definition of quasirandomness: a regular or approximately regular graph is 
called quasirandom if it contains approximately the expected number of 
$C_4$s that would be in a random graph chosen with the edge probability 
dictated by the density of the graph (see for example \cite{CGW89}, 
\cite{T87}). The standard 
results of the theory imply that if a graph satisfies this criterion then 
it also satisfies many of the properties that are expected with high 
probability of a random graph. For example, and this is what will be 
important to us, it contains approximately the expected number of any 
small graph.

The properties that we now know about a colouring of a $K_n$ not 
containing either a red $K_{k+1}$ or a blue $K_{l+1}$ are enough to tell 
us that both the red and blue components of our colouring are 
quasirandom, so we see that in such a colouring we have approximately the 
expected number of any small graph in either colour. In particular, for 
any fixed $r$, we have approximately 
\[\left( \frac{k}{k+l} \right)^{\binom{r}{2}} n^r\] 
ordered red $r$-tuples (we find it more convenient to count $r$-tuples 
rather than $K_r$s since we then don't have to worry about multiple 
counting in our estimates, but it might perhaps be best to think of it in 
terms of counting $K_r$s). 

If this were in fact precise then it would be inconsistent with the fact 
that any red $(r-1)$-tuple lies in at most $r(k - r + 2, l+1) - 
1$ red $r$-tuples, since this gives an upper bound on the number of red 
$r$-tuples of
\[r(k - r + 2, l+1) \left( \frac{k}{k+l} \right)^{\binom{r-1}{2}} 
n^{r-1},\]
which, since 
\[r(k - r + 2, l+1) \leq \frac{f(k-r+1,l)}{f^*(k,l)} 
\frac{k\cdots(k-r+2)}{(k+l)\cdots(k+l-r+2)} n,\]
is strictly less than the expected number if the rate of change of 
$f$ is sufficiently small.

It is precisely this contradiction which allows us to prove our result. 
There are of course several technical caveats, the most interesting of 
which is that, in order to derive Theorem 1.1, it is not sufficient to 
know that the graph is simply quasirandom. It is necessary to apply instead 
our local condition that we have approximately the expected number of red 
$C_4$s across any given edge. Theorem 1.2, on the other hand, is derivable 
from the quasirandomness condition alone (though we will not do this here). 

Secondly, the argument as stated above is slightly illusory - in order to 
derive a useful result it is necessary to take into account the fact that 
a change in the number of $K_{r-1}$s will be reflected by a change in the 
number of $K_r$s. Without doing this, we would be able to do no better 
than Thomason's result. 

Where, incidentally, do we depart from Thomason's work? His proof is 
essentially the argument given above in the case $r=3$. He counts, in two 
different ways, the number of monochromatic triangles within a graph not 
containing a red $K_{k+1}$ or a blue $K_{l+1}$, showing that, unless a 
bound of the form 
\[r(k+1, l+1) \leq \exp \left\{ -\frac{l}{2k} \mbox{log} k + A 
\sqrt{\log k} \right\} \binom{k+l}{k}\]
held, there would be a contradiction. While his method of finding an 
upper bound for the number of monochromatic triangles is similar to ours 
above (the number of red triangles across a given red edge is at most 
$r(k-1, l+1) - 1$, and we know, approximately, the number of red edges), 
his method for finding a lower bound is to apply Goodman's formula 
\[T = \frac{1}{2} \left[ \sum_x \binom{d_x}{2} + \sum_x 
\binom{n-1-d_x}{2} - \binom{n}{3} \right],\]
where by $d_x$ we mean the red degree of the vertex $x$. This formula is 
only dependent upon the degree sequence, and so, knowing that every 
degree is approximately what's expected, we can show that the number of 
monochromatic triangles is approximately what's expected. Our main 
advance then is to have shown how we can use the quasirandomness 
conditions to circumvent the fact that there is no Goodman-type formula 
for $r \geq 4$. 

This discussion raises one further question: why, if Thomason's result 
implies an off-diagonal estimate as well as a diagonal one, and our 
arguments are the natural extension of Thomason's argument, do we not 
also have off-diagonal theorems which include Theorem 1.1 and Theorem 1.2 
as special cases? The first part of the answer is that for Theorem 1.2 we 
do, the following theorem being our main result in this case:

\begin{theorem} 
Let $s$ and $\epsilon$ be fixed positive constants with $\epsilon \leq 
1$. Then there exists a constant $C_{s, \epsilon}$ such that for $k \geq 
l \geq \epsilon k$ 
\[r(k+1, l+1) \leq C_{s, \epsilon} \exp\left\{-s \frac{l}{k} \log 
k\right\} \binom{k+l}{k}.\]
\end{theorem}

However, if we now fix $\epsilon$ and let $s$ increase, the best theorem 
that we can deduce from our knowledge of the growth rate of the $C_{s, 
\epsilon}$ is

\begin{theorem} 
For all $\epsilon \leq 1$ there exists a constant $C_\epsilon$ such that 
for all $k$ and $l$ with $k \geq l \geq \epsilon k$, 
\[r(k+1, l+1) \leq \exp\{- C_\epsilon \log^{3/2} k \} \binom{k+l}{k}.\]
\end{theorem}

Interestingly, as with Theorem 1.2, we only need the ordinary 
quasirandomness condition and not its local counterpart to derive these 
results. It is only when we make $\epsilon$ tend towards $1$ as $s$ 
increases that our local conditions become genuinely useful. So the 
reason why we don't have an off-diagonal version of Theorem 1.1 now 
emerges: our method doesn't allow one. 

We will not prove Theorems 2.1 and 2.2 in this paper, concentrating 
instead on the diagonal results. It should however be clear to the reader 
how we can go about changing our main results in order to derive them.

We begin the proof proper in the next section by considering, more 
formally, the various regularity conditions that a graph containing 
neither a red $K_{k+1}$ nor a blue $K_{l+1}$ must satisfy, and showing 
what these conditions imply about such a graph.

\section{The Regularity Conditions} 

The following notation will prove essential to us in what follows:

\demo{Definition} Suppose we have a red/blue colouring of the edges 
of the complete graph on $n$ vertices, and let $V$ be the set of 
vertices. Then we define the balanced function of the colouring around 
probability $p$ as the function $g: V \times V \rightarrow \mathbb{R}$ 
with 
\[g(x, y) = A(x,y) - p,\]
where $A: V \times V \rightarrow \mathbb{R}$ is the characteristic 
function of red edges, that is $A(x,y)$ is $1$ if there is a red edge 
between $x$ and $y$ and $0$ otherwise.
\Enddemo

Note that normally one chooses the probability $p$ in such a 
way as to make $\sum_{x,y} g(x,y) = 0$, but for the sake of simplicity in 
our exposition, we will be centring around a probability which is not quite 
the correct balanced probability, but which is very close.

We will also need to introduce two constants, $\gamma$ and $\delta$, 
which bound the growth (or rather fall) of $f(k,l)$ with respect to $k$ 
and $l$ respectively. Our main result in the next section will be an 
inequality telling us what kind of rate of change of $f(k,l)$ is 
admissible. More specifically, we will assume that we have two real 
numbers $\gamma$ and $\delta$ and a natural number $n = \lfloor f(k,l) 
\binom{k+l}{k} \rfloor = f^*(k,l) \binom{k+l}{k}$, such that for $m = 1, 
2$ and $r-1$, each of the inequalities
\[r(k+1 -m, l+1) \leq f(k-m, l) \binom{k-m+l}{k-m},\]
\[r(k+1, l+1-m) \leq f(k, l-m) \binom{k+l-m}{k},\]
\[\frac{f(k-m,l)}{f^*(k,l)} \leq 1 + m \gamma \mbox{ and }
\frac{f(k,l-m)}{f^*(k,l)} \leq 1 + m \delta\]
holds. What we will show (by the counting $K_r$s argument we discussed in 
the last section) is that if $k \geq l$, where $k$ and $l$ are 
sufficiently large numbers of approximately the same magnitude, and if 
\[k\gamma + l\delta \leq \frac{r-3}{2} \frac{l}{k},\]
then 
\[r(k+1, l+1) \leq f(k,l) \binom{k+l}{k}.\]
The conditions on $\gamma$ and $\delta$ essentially amount to $\gamma$ 
and $\delta$ being the partial derivatives of $\phi(k,l) = -\log f(k,l)$ 
with respect to $k$ and $l$ respectively. Thus, if we consider the 
inequality $k\gamma + l\delta \leq \frac{r-3}{2} \frac{l}{k}$ as a 
partial differential equation (by putting $\gamma = \frac{\partial 
\phi}{\partial k}$ and $\delta = \frac{\partial \phi}{\partial l}$), it 
is easy to see that taking $f(k,l) = \exp\{-\frac{r-3}{2} \frac{l}{k} 
\log k\}$ for $k \geq l$ works as a potential solution. Indeed a more 
careful treatment of this argument, taking into account the fact that 
$\gamma$ and $\delta$ do not quite equal the respective derivatives, is 
what will allow us to derive our results.

The specifics of this must, however, wait until later sections. The task 
at hand is show what we can say about large graphs not containing either red 
$K_{k+1}$s or blue $K_{l+1}$s. We begin by writing our various 
regularity conditions as constraints on the size of certain products of 
the balanced function:

\begin{lemma}
Let $k$ and $l$ be natural numbers, let $\gamma$ and 
$\delta$ be real numbers and let $n = \lfloor f(k,l) 
\binom{k+l}{k} \rfloor = f^*(k,l) \binom{k+l}{k}$. Suppose that for $m = 
1$ and $m = 2$ each of the inequalities
\[r(k+1 -m, l+1) \leq f(k-m, l) \binom{k-m+l}{k-m},\]
\[r(k+1, l+1-m) \leq f(k, l-m) \binom{k+l-m}{k},\]
\[\frac{f(k-m,l)}{f^*(k,l)} \leq 1 + m \gamma \mbox{ and }
\frac{f(k,l-m)}{f^*(k,l)} \leq 1 + m \delta\]
holds.

Then, in any red/blue colouring of $K_n$ not containing either a red 
$K_{k+1}$ or a blue $K_{l+1}$, the balanced function $g(x,y)$ of the 
colouring around $p = \frac{k}{k+l}$ satisfies 
\[-\frac{l \delta}{k+l} n \leq \sum_y g(x,y) \leq \frac{k \gamma}{k+l} n\]
for all $x$, and
\[\sum_y g(x,y)g(y,z) \leq 2 \frac{\max(k,l)}{(k+l)^2} (k\gamma + 
l\delta) n + 1\]
for all $x$ and $z$ with $x \neq z$. 
\end{lemma}

\Proof
The first part of the lemma follows from our observation in Section 2 
that for any vertex $x$ in our colouring we have
\[\left( 1 - \frac{f(k, l-1)}{f^*(k,l)} \frac{l}{k+l} \right) n\leq d_x 
\leq \frac{f(k-1, l)}{f^*(k,l)} \frac{k}{k+l} n.\]
Noting that $d_x = \sum_y A(x,y)$, $A(x,y) = p + g(x,y)$, and applying 
our assumptions on the growth rate of $f$ gives the required result. To 
prove the upper bound, for example, note that
\begin{eqnarray*}
\frac{k}{k+l} n + \sum_y g(x,y) & = & d_x \\
& \leq & \frac{f(k-1, l)}{f^*(k,l)} \frac{k}{k+l} n \\
& \leq & (1 + \gamma) \frac{k}{k+l} n. 
\end{eqnarray*}
Subtracting $\frac{k}{k+l} n$ from either side then gives the required 
bound.

For the second part of the lemma note that no red edge $(x, z)$ can lie in 
more than $r(k-1, l+1) - 1$ red triangles. This implies that
\[\sum_y  A(x, y)A(y, z) \leq r(k-1, l+1) - 1.\]
If we split up the left-hand side we then get, by using the conditions of 
the theorem, that
\[p^2 n + p \sum_y g(x,y) + p \sum_y g(y, z) + \sum_y g(x, y)g(y, z) \leq 
p^2 (1 + 2 \gamma) n,\]
and hence by the first part of the lemma
\[\sum_y g(x, y)g(y, z) \leq 2 \frac{k}{(k+l)^2} (k\gamma + l\delta) n.\]
The result follows similarly for blue edges, although we need to be a 
little bit careful, since we get two extra degenerate ``triangles" (those 
for which $y = x$ or $y=z$). 
\Endproof

In counting the number of red $K_r$s in a given colouring, we will use 
the following notation: 

\demo{Notation} Fix a red/blue colouring on $K_n$ and let $g(x,y)$ be the 
balanced function of the colouring around probability $p$. Suppose also 
that $K_r$ is the complete graph on the $r$ vertices $v_1, v_2, \cdots, 
v_r$, with $r \leq n$. Then, for every subgraph $H$ of $K_r$,we write 
\[g_H = \sum_{x_1, \cdots, x_r} \prod_{(v_i, v_j) \in E(H)} g(x_i, x_j),\]
where the sum is taken over all $r$-tuples of vertices in $K_n$ 
(including degenerate terms where two or more of the $x_i$ are the same).
\Enddemo

By rights this is a function of $n$ and $r$ as well as $H$, but we will 
be almost universally consistent about counting $K_r$s within $K_n$s, so 
these labels are essentially redundant.

Given this notation, the number of red $K_r$s (or rather red $r$-tuples) 
in a colouring of $K_n$ is given by
\begin{eqnarray*}
\sum_{x_1, \cdots, x_r} \prod_{(v_i, v_j) \in E(K_r)} A(x_i, x_j) & = &
\sum_{x_1, \cdots, x_r} \prod_{(v_i, v_j) \in E(K_r)} (p + g(x_i, x_j)) \\
& = & \sum_{H \subset K_r} p^{\binom{r}{2} - e(H)} g_H,
\end{eqnarray*}
where, again, the sum is taken over all $r$-tuples of vertices in $K_n$. 
So in order to estimate the number of $K_r$s we will need to be able to 
estimate $g_H$ for each and every subgraph $H$ of $K_r$. Almost all of 
the estimates we will need are encapsulated in the next lemma, which 
shows how we may use our local quasirandomness condition to obtain 
estimates on products of the balanced function.

Utilising the information provided by the previous lemma, we shall now 
assume that we have $\sum_y g(x,y)g(y,z) \leq \nu n$ for all $x$ and $z$ 
with $x \neq z$, where $\nu$ is some positive constant. The next 
lemma tells us that if $H$ has a vertex of degree $d$ then (to the 
highest order in $n$) $|g_H| \leq \sqrt{2} \nu^{d/2} n^r$. 
Within the statement of the lemma, we will make the simple assumption 
that $\nu \leq 1$. This is not strictly necessary but tidies up the form 
of the lemma, and as we shall see later is trivially satisfied for $k$ 
and $l$ large.

\begin{lemma}
Suppose that the balanced function $g(x,y)$ of a red/blue colouring of a 
graph on $n$ vertices satisfies
\[\sum_y g(x, y)g(y, z) \leq \nu n\]
for all $x$ and $z$ with $x \neq z$, and some fixed positive real $\nu$. 
Then, provided that $\nu \leq 1$,
\[|\sum_y \sum_{x_1, \cdots, x_{c+d}} g(y, x_1) \cdots g(y, x_d) h(x_1, 
\cdots, x_{c+d})| \leq \sqrt{2} \nu^{d/2} n^{c+d+1} + \frac{1}{\sqrt{2} 
\nu^{d/2 + 1}} n^{c+d},\]
for any function $h$ of $c+d$ vertices which is bounded above in absolute 
value by 1.
\end{lemma}

\Proof
For $d$ odd, we have
\begin{eqnarray*}
&&|\sum_y \sum_{x_1, \cdots, x_{c+d}} g(y, x_1) g(y, x_2) \cdots g(y, x_d) 
h(x_1, \cdots, x_{c+d})|^2\\
& \leq & n^{c+d} \sum_{x_1, \cdots, x_{c+d}} | \sum_y 
g(y, x_1) g(y, x_2) \cdots g(y, x_d) h(x_1, \cdots, x_{c+d})|^2 \\
& \leq & n^{c+d} \sum_{x_1, \cdots, x_{c+d}} | \sum_y g(y, x_1) g(y, x_2) 
\cdots g(y, x_d)|^2 \\
& = & n^{2c+d} \sum_{y, y'} (\sum_x g(y, x) g(x, y'))^d \\
& \leq & \nu^d n^{2c + 2d + 2} + n^{2c+2d+1}, 
\end{eqnarray*}
where the remainder comes from the degenerate terms. Since this is less 
than the square of 
\[\nu^{d/2} n^{c+d+1} + \frac{1}{2 \nu^{d/2}} n^{c+d},\]
we are done in this case.

For $d$ even, the proof is the same until we reach the second last line, 
when we need to estimate
\[\sum_{y, y'} (\sum_x g(y, x) g(x, y'))^d.\]
To do this we split our sum into two pieces, a set $P$ of 
edges $(y, y')$
where $\sum_x g(y, x) g(x, y')$ is positive and a similar set $N$ where 
$\sum_x g(y, x) g(x, y')$ is negative. Then the proof in the odd case 
tells us, since a sum of squares is positive, that
\[\sum_{(y, y') \in P} (\sum_x g(y, x) g(x, y'))^{d+1} \geq -\sum_{(y, y') 
\in N} (\sum_x g(y, x) g(x, y'))^{d+1} \]
which implies
\begin{eqnarray*}
\sum_{y, y'} |\sum_x g(y, x) g(x, y')|^{d+1} & \leq & 2 \sum_{(y, y') \in 
P} (\sum_x g(y, x) g(x, y'))^{d+1} \\
& \leq & 2 \nu^{d+1} n^{d+3} + 2 n^{d+2}.
\end{eqnarray*}
Finally, applying the power mean inequality, we get
\begin{eqnarray*}
\sum_{y, y'} (\sum_x g(y, x) g(x, y'))^d & \leq & n^{\frac{2}{d+1}}  
(\sum_{y, y'} |\sum_x g(y, x) g(x, y')|^{d+1})^{\frac{d}{d+1}} \\
& \leq & n^{\frac{2}{d+1}} (2 \nu^{d+1} n^{d+3} + 2 n^{d+2})^{\frac{d}{d+1}} \\
& \leq & 2 \nu^{d} n^{d+2} + \frac{2}{\nu} n^{d+1},
\end{eqnarray*}
so we are done in this case as well.
\Endproof

The result mentioned before the lemma now follows from taking $y$ to be a 
vertex within $H$ of degree $d$. The function $h$ is then what remains, 
i.e. a certain product of balanced functions, and so satisfies the 
requirement of the lemma. 

Ultimately, as we shall see in the next section, we would like to show 
that as many $g_H$ terms as possible vanish to more than the first order 
in $\gamma$ and $\delta$. While the above results are sufficient to show 
that this is so when the graph $H$ has maximum degree $3$ or more, it 
still leaves a large collection of graphs of maximum degree $2$ for which 
we have not reached this bound. The next lemma shows, however, that if we 
use the degree-regularity condition as well as the quasirandomness 
condition, then we have the required bounds except in the two cases where 
$H$ is a $K_2$ or a $K_3$.

\begin{lemma}
Suppose that the balanced function $g(x,y)$ of a red/blue colouring of a 
graph on $n$ vertices satisfies
\[|\sum_y g(x, y)| \leq \mu n\]
for all $x$, and
\[\sum_y g(x, y)g(y, z) \leq \nu n\]
for all $x$ and $z$ with $x \neq z$, and some fixed positive constants $\mu$
and $\nu$ with $\nu \leq 1$. Then, for $l \geq 3$, 
\[|\sum_{x_1, \cdots, x_l} g(x_1, x_2) g(x_2, x_3) \cdots g(x_{l-1}, x_l)| 
\leq 2 \mu^{l+1 - 2 \lfloor l/2 \rfloor} \nu^{\lfloor l/2 \rfloor - 1} 
n^l + \frac{2 \mu^{l+1 - 2 \lfloor l/2 \rfloor}}{\nu^3} n^{l-1} \]
and
\[|\sum_{y_1, \cdots, y_l} g(y_1, y_2) g(y_2, y_3) \cdots g(y_l, y_1)| \leq 
2 \nu^{\lfloor l/2 \rfloor} n^l + \frac{2}{\nu} n^{l-1}.\]
\end{lemma}

\Proof
For the first part we simply apply the H\"{o}lder inequality:
\begin{eqnarray*} 
&&|\sum_{x_1, \cdots, x_l} g(x_1, x_2) \cdots g(x_{l-1}, x_l)|^{\lceil l/2 
\rceil}\\
& = &  \left| \sum_{x_2, x_4, \cdots} \left(\sum_{x_1} g(x_1, x_2) \right) 
\left(\sum_{x_3} g(x_2, x_3)g(x_3, x_4) \right) \cdots \right|^{\lceil l/2 
\rceil} \\
& \leq & \left(\sum_{x_2, x_4, \cdots} |\sum_{x_1} g(x_1, x_2)|^{\lceil
l/2 \rceil} \right) \left(\sum_{x_2, x_4, \cdots} |\sum_{x_3} g(x_2,
x_3)g(x_3, x_4)|^{\lceil l/2 \rceil} \right) \cdots \\
& \leq & \left( \mu^{\lceil l/2 \rceil} n^l\right)^{l+1 - 2 \lfloor l/2 
\rfloor} \left( 2 \nu^{\lceil l/2 \rceil} n^l + \frac{2}{\nu} n^{l-1} 
\right)^{\lfloor l/2 \rfloor - 1} \\
& \leq & \left( 2 \mu^{l+1 - 2 \lfloor l/2 \rfloor} \nu^{\lfloor l/2 
\rfloor - 1} n^l \right)^{\lceil l/2 \rceil} \left(1 + 
\frac{1}{\nu^{\lceil l/2 \rceil + 1} n} \right)^{\lceil l/2 \rceil},
\end{eqnarray*} 
which implies the result. The second part follows similarly.
\Endproof

\section{The Fundamental Lemma}

In this section we will prove an extension of a lemma due to Thomason 
\cite{T88} which gives an inequality telling us how quickly our function 
$f(k,l)$ may change. The main idea of our proof is one that we have 
already seen. Instead of counting the number of monochromatic triangles 
as Thomason did, we will count the number of monochromatic $K_r$s (or 
rather a certain weighted sum of the number of red $K_r$s and the number 
of blue $K_r$s), showing, using the fact that our graph must be 
random-like if it does not contain the required cliques, that this is 
approximately what is expected. On the other hand, we can again bound the 
number of monochromatic $K_r$s above using the following further 
generalisation of the Erd\H{o}s-Szekeres condition: in a graph not 
containing a red $K_{k+1}$ or a blue $K_{l+1}$ any red $K_{r-1}$ is 
contained in at most $r(k-r+2, l+1) - 1$ red $K_r$s and any 
blue $K_{r-1}$ is contained in at most $r(k+1, l-r+2) - 1$ blue 
$K_r$s. Then since the number of $K_{r-1}$s can also be estimated (as 
approximately the expected number) we have an upper bound which we can 
balance against our lower bound. 

Again, as we mentioned in the outline, it will be necessary in the proof 
to take into account the fact that the number of $K_r$s and the number of 
$K_{r-1}$s are not independent of one another, being composed almost 
entirely of like terms, although in different proportions. While most of 
these terms may be reduced to $o(1)$ factors at the outset as being quite 
unimportant to the argument, the terms coming from single edges and 
triangles, which are the highest order, and hence the critical, terms, 
will be left unestimated until after we have balanced the number of red 
$K_r$s against $r(k-r+2,l+1)$ times the number of red $K_{r-1}$s. Doing 
this allows us to reduce the error term coming from the single edges from 
being of the order of $r^2 \sum_{x,y} g(x,y)$ to being $r \sum_{x,y} 
g(x,y)$, since the single edge terms which occur in counting the number 
of $K_{r-1}$s cancel out most of the like terms which we get in counting 
the number of $K_r$s. Without this care, our result would yield no 
improvement over the old bound. 

Before we begin, we need to present a few more remarks in order to 
illuminate some of the assumptions of the lemma. 
What we will prove is that if $k$ and $l$ are sufficiently large 
depending on $r$, $k \geq l \geq \left(1 - \frac{1}{r}\right) k$ and 
\[k \gamma + l \delta \leq \frac{r-3}{2} \frac{l}{k},\]
then (given the obvious induction hypothesis), we have 
\[r(k+1, l+1) \leq f(k,l) \binom{k+l}{k}.\]

Now, as at the start of Section 3, we see that with an inequality of this 
form, we expect $f(k,l)$ to be roughly of the form
\[\exp\left\{- \frac{r-3}{2} \frac{l}{k} \log k\right\},\]
or some multiple thereof. One result of this is that we expect both 
$|\gamma|$ and $|\delta|$ to be bounded by $\frac{r-3}{2} \frac{\log 
k}{k}$. Since our eventual hope is to prove that $f(k,l)$ has such a form 
we will in the course of our forthcoming proof, in order to simplify the 
final form of the result, make the assumptions that $f(k,l)$ is at the 
smallest equal to $\exp \{-r \frac{l}{k} \log k \}$ and that both 
$|\gamma|$ and $|\delta|$ are smaller than $r \frac{\log k}{k}$. There is 
no deep mystery to our using $r$ rather than $\frac{r-3}{2}$ here. It's 
just neater, and makes the lemma look slightly more digestible.

We are now ready to begin the formalities.

\begin{lemma}
Let $r$ be a natural number, let $\gamma$ and $\delta$ be real numbers 
and let $n = \lfloor f(k,l)\binom{k+l}{k} \rfloor = f^*(k,l) 
\binom{k+l}{k}$. Suppose that, for 
$m=1$, $m=2$ and $m=r-1$, each of the inequalities
\[r(k+1 -m, l+1) \leq f(k-m, l) \binom{k-m+l}{k-m},\]
\[r(k+1, l+1-m) \leq f(k, l-m) \binom{k+l-m}{k},\]
\[\frac{f(k-m,l)}{f^*(k,l)} \leq 1 + m \gamma \mbox{ and }
\frac{f(k,l-m)}{f^*(k,l)} \leq 1 + m \delta\]
holds.
Suppose also that 
\begin{description}
\item{1.}
$k \geq l \geq \left( 1 - \frac{1}{r} \right) k$,
\item{2.}
$|\gamma|$ and $|\delta|$ are both smaller than $r 
\frac{\log k}{k}$ and 
\item{3.}
$f(k,l) \geq \exp \{-r \frac{l}{k} \log k \}$.
\end{description}
Then there exists a constant $c$ such that if $k$ and $l$ are both 
greater than $r^{cr}$, and 
\[k \gamma + l \delta \leq \frac{r-3}{2} \frac{l}{k},\]
the inequality
\[r(k+1, l+1) \leq n \leq f(k,l) \binom{k+l}{k}\]
holds.
\end{lemma}

\Proof
To begin, note that, from Lemma 3.1, in a colouring avoiding red 
$K_{k+1}$s and blue $K_{l+1}$s, we must have that the balanced function 
$g(x,y)$ satisfies
\[-\frac{l \delta}{k+l} n \leq \sum_y g(x,y) \leq \frac{k \gamma}{k+l} n\]
for all $x$, and therefore, using assumption $2$ of the lemma, we have that
\[|\sum_y g(x,y)| \leq r \frac{\log k}{k} n.\]
Also from Lemma 1, note that, since $k\gamma + l \delta \leq 
\frac{r-3}{2}$, we have that 
\[\sum_y g(x,y) g(y,z) \leq \frac{r-3}{k+l} n\]
for all $x$ and $z$ with $x \neq z$ (we may subsume the O(1) term into 
the $n$ term for $k$ and $l$ larger than some fixed constant - it is in 
performing this kind of estimate that we will use assumption 3 of the lemma).

For later brevity we will use the notation that $\sum_{x,y} g(x,y) = 
\frac{s}{k+l} n^2$, and we will also write $\sum_{x,y,z} g(x,y) g(y,z) 
g(z,x) = \frac{t}{k+l} n^3$.
Moreover, we will denote the quantity $\frac{r-3}{k+l}$ by $\nu$ so that 
\[\sum_y g(x,y) g(y,z) \leq \nu n,\]
noting that, for $k+l \geq r$, we have $\nu \leq 1$.

Recall that the number of $r$-tuples spanning a red clique is given by
\[\sum_{H \subset K_r} p^{\binom{r}{2} - e(H)} g_H.\]
Our first aim will be to show that the contribution of all terms in this 
sum other than the main term (corresponding to the null set), the edge 
terms and the triangle terms can be made smaller in absolute value than 
$\frac{1}{2^{\binom{r}{2}} r^{d r} k}$ for any fixed $d$, by taking $k$ 
and $l$ to be larger than $r^{cr}$ for some appropriately large $c$.

Let us denote by $S$ the set of subgraphs of $K_r$ other than the null 
graph, the edges, and the triangles. We will split this set into two 
further subsets, $S'$, the set of all subgraphs with maximum degree 
greater than or equal to $3$, and $S''$, the complement of this set in 
$S$. 

For graphs in $S'$, Lemma 3.2 tells us that, for $k$ and $l$ greater than 
$r^{cr}$,
\begin{eqnarray*}
|g_H| & \leq & \sqrt{2} \nu^{\Delta/2} n^r + \frac{1}{\sqrt{2} 
\nu^{\Delta/2 + 1}} n^{r-1} \\ 
& \leq & \sqrt{2} \left( \frac{r}{k} \right)^{\Delta/2} n^r + 
\frac{1}{\sqrt{2}} \left(\frac{k}{r}\right)^{\Delta/2 + 1} n^{r-1} \\
& \leq & \frac{1}{r^{c_1 \Delta r} k} n^r,
\end{eqnarray*} 
where $\Delta$ is the maximum degree of $H$, and where $c_1$ depends on 
and grows with $c$. 

Now, every graph in $S''$ either contains a path of length two or a cycle 
of length 4, in which case we have, from Lemma 3.3, and using our bounds 
on $\mu = \max_x |\sum_y g(x,y)|$ and $\nu$, that
\begin{eqnarray*}
|g_H| & \leq & 2 r^2 \frac{\log^2 k}{k^2} n^r + 2 \left( \frac{k}{r} 
\right)^3 n^{r-1} \\
& \leq & \frac{1}{r^{c_2 r} k} n^r,
\end{eqnarray*}
or is a product of single edges and triangles, in which case
\begin{eqnarray*}
|g_H| & \leq & 4 r^2 \frac{\log^2 k}{k^2} n^r + 8 n^{r-1} + 4 \left( 
\frac{k}{r} \right)^2 n^{r-2} \\
& \leq & \frac{1}{r^{c_2 r} k} n^r,
\end{eqnarray*}
where again $c_2$ is just some constant that grows with $c$. 

Before we proceed with our estimate we also need to note firstly that the 
number of graphs with maximum degree $\Delta$ or less is at most 
$r^{\Delta r}$ and also that the maximum number of edges in such a graph 
is $\Delta r$ (we may of course divide by a 2 here but this is not 
necessary for our estimates). 

We now have, using the fact that $p$ is between $1/3$ and $1/2$, that
\begin{eqnarray*}
\sum_{H \in S} p^{\binom{r}{2} - e(H)} |g_H| & \leq & \left( \frac{1}{2} 
\right)^{\binom{r}{2}} \sum_{H \in S'} \frac{3^{e(H)}}{r^{c_1 \Delta r} 
k} n^r + \sum_{H \in S''} \frac{3^{e(H)}}{r^{c_2 r} k} n^r \\
& \leq & \left( \frac{1}{2} \right)^{\binom{r}{2}} \sum_{\Delta=3}^r 
\frac{(3r)^{\Delta r}}{r^{c_1 \Delta r} k} n^r + \frac{(3r)^{2r}}{r^{c_2 
r} k} n^r \\
& \leq & \frac{1}{2^{\binom{r}{2}} r^{d r} k} n^r
\end{eqnarray*}
for $c$ chosen sufficiently large depending on $d$.

So, getting back to our original intentions, we see that the number of 
$r$-tuples spanning a red $K_r$ is greater than or equal to
\[\frac{k^{\binom{r}{2}}}{(k+l)^{\binom{r}{2}}} n^r + \binom{r}{2} 
\frac{k^{\binom{r}{2} - 1}}{(k+l)^{\binom{r}{2}}} s n^{r}
+ \binom{r}{3} \frac{k^{\binom{r}{2}-3}}{(k+l)^{\binom{r}{2}-2}} t n^{r} 
- \frac{1}{2^{\binom{r}{2}} r^{dr} k} n^r.\]

On the other hand we have that the number of $r$-tuples with a red $K_r$ 
across it is less than the number of $(r-1)$-tuples with a red $K_{r-1}$ 
across it times $r(k+2-r, l+1)$. So we have that the number of $K_r$s is 
at most 
\[\left(\frac{k^{\binom{r-1}{2}}}{(k+l)^{\binom{r-1}{2}}} n^{r-1} + 
\binom{r-1}{2} \frac{k^{\binom{r-1}{2} - 1}}{(k+l)^{\binom{r-1}{2}}} 
s n^{r-1} + \binom{r-1}{3} \frac{k^{\binom{r-1}{2} - 
3}}{(k+l)^{\binom{r-1}{2} - 2}} t n^{r-1}\right.\] 
\[\left.+ \frac{1}{(r-1)^{d(r-1)} k} n^{r-1}\right)
\times \left( \frac{k(k-1)\cdots(k-r+2)}{(k+l)(k+l-1)\cdots(k+l-r+2)}\right)
(1+(r-1)\gamma) n.\]

Now we are going to subtract the lower bound from the upper bound, and
divide through by $n^r$ to get an
inequality which must hold if the graph on $n$ vertices contains neither a
red $K_{k+1}$ or a blue $K_{l+1}$. In so doing it is necessary to use the 
fact that
\begin{eqnarray*}
&& \left( \frac{k}{k+l} \right)^{r-2} - \left(\frac{(k-1)\cdots(k-r+2)}{(k+l-1)
\cdots(k+l-r+2)}\right)\\
& \geq & \binom{r-1}{2} \frac{l k^{r-3}}{(k+l)(k+l-1)\cdots(k+l-r+2)} 
- \frac{2^r}{(k+l)^2}\\
& \geq & \binom{r-1}{2} \frac{l k^{r-3}}{(k+l)^{r-1}} 
- \frac{2^r}{(k+l)^2}
\end{eqnarray*} 
so that all of the second order terms arising from the use of 
$\left(\frac{(k-1)\cdots(k-r+2)}{(k+l-1)\cdots(k+l-r+2)}\right)$ instead 
of $\left( \frac{k}{k+l} \right)^{r-2}$, other than that coming from the 
main term, will in fact be $1/k^2$ terms or smaller still. They can 
therefore be subsumed into the remainder term. This yields
\[\binom{r-1}{2} 
\frac{lk^{\binom{r}{2}-1}}{(k+l)^{\binom{r}{2}+1}}
+ (r-1) \frac{k^{\binom{r}{2} - 1}}{(k+l)^{\binom{r}{2}}} s + 
\binom{r-1}{2} \frac{k^{\binom{r}{2}-3}}{(k+l)^{\binom{r}{2}-2}}t\] 
\[-(r-1)\frac{k^{\binom{r}{2}}\gamma}{(k+l)^{\binom{r}{2}}} \leq 
\frac{1}{2^{\binom{r}{2}} r^{Dr} k},\]
where $D$ is again just some constant which grows with $c$. 

Similarly, if we count blue $K_r$s, though we have to be a little bit 
careful about degenerate terms and the fact that $1-p$ may be slightly 
bigger than $1/2$, we get 
\[\binom{r-1}{2} \frac{kl^{\binom{r}{2}-1}}{(k+l)^{\binom{r}{2}+1}}
- (r-1) \frac{l^{\binom{r}{2} - 1}}{(k+l)^{\binom{r}{2}}} s - 
\binom{r-1}{2} \frac{l^{\binom{r}{2}-3}}{(k+l)^{\binom{r}{2}-2}}t\]
\[-(r-1)\frac{l^{\binom{r}{2}}\delta}{(k+l)^{\binom{r}{2}}} \leq 
\frac{1}{2^{\binom{r}{2}} r^{Dr} k}.\]

Now we take the weighted sum of these inequalities in such a way as to 
make the triangle terms disappear, by adding $k^{\binom{r}{2} - 3}$ times 
the first inequality to $l^{\binom{r}{2} - 3}$ times the second, to get
\[\binom{r-1}{2} \frac{(kl)^{\binom{r}{2}-2}}{(k+l)^{\binom{r}{2}}} + 
(r-1)\frac{(k^2 - l^2)(kl)^{\binom{r}{2}-3}}{(k+l)^{\binom{r}{2}}}s -(r-1) 
\frac{(k^3 \gamma + l^3 
\delta)(kl)^{\binom{r}{2}-3}}{(k+l)^{\binom{r}{2}}} \leq \frac{2 
k^{\binom{r}{2} - 4}}{2^{\binom{r}{2}} r^{Dr}}.\]
Finally, provided that $l \geq \left(1 - \frac{1}{r} \right) k$ and $c$ has been chosen large enough, we see that we can subsume the error term on the right hand side into the first term on the left hand side, the first term on the left hand side being then larger than
\[\binom{r-1}{2} \frac{k^{\binom{r}{2} - 4}}{e^{2r} 2^{\binom{r}{2}}}\]
(here we have used the fact that $1 - \frac{1}{r} \geq e^{-2/r}$ for $r 
\geq 2$).

Therefore, subsuming this term, we see that
\[\frac{(r-1)(r-3)}{2} \frac{(kl)^{\binom{r}{2}-2}}{(k+l)^{\binom{r}{2}}} 
+ (r-1)\frac{(k^2 - l^2)(kl)^{\binom{r}{2}-3}}{(k+l)^{\binom{r}{2}}}s - 
(r-1) \frac{(k^3 \gamma + l^3   
\delta)(kl)^{\binom{r}{2}-3}}{(k+l)^{\binom{r}{2}}} < 0.\]
Simplifying gives
\[\frac{(r-3)}{2} kl + (k^2 - l^2) s - (k^3 \gamma + l^3 \delta) < 0.\] 
Recall now that
\[s = \frac{(k+l) \sum_{x,y} g(x,y)}{n^2} \geq -l \delta,\]
so therefore, since $k \geq l$,
\[k \gamma + l \delta > \frac{r-3}{2} \frac{l}{k}.\]
This contradicts the assumptions of the lemma, and so we are done.
\Endproof

\section{Using the Inequality}

All that now remains to be done is to find a function that satisfies the 
conditions of Lemma 4.1. The basic idea is to note that if we choose a 
continuously differentiable function $\alpha: [0, \infty) \rightarrow [0, 
\infty)$, then the function
\[f(k,l) = \exp\left\{ - \alpha(l/k) \log (k+l) \right\}\]
satisfies the equation
\[k \gamma' + l \delta' = \alpha(l/k),\]
where by $\gamma'$ and $\delta'$ we mean the derivatives of $-\log 
f(k,l)$ with respect to $k$ and $l$.  

To use this fact we will choose a function $\alpha_r$ which is everywhere 
less than or equal to the function $\kappa_r$, where
\[\kappa_r (x) = \left\{ \begin{array}{ll}
0 & \mbox{ if $0 \leq x < 1 - \frac{1}{r}$;}\\
\frac{r-3}{2} x & \mbox{ if $1 - \frac{1}{r} \leq x \leq 1$;}\\
\kappa_r(1/x) & \mbox{ if $x \geq 1$.} \end{array} \right.\]
and which, moreover, is twice-differentiable. The specific function, if 
it chosen appropriately, will then be such that the true $\gamma$ and 
$\delta$ differ by very little from $\gamma'$ and $\delta'$ for $k$ and 
$l$ chosen quite large, and this will allow us to conclude, for a 
suitably chosen $\alpha_r$, that
\[k \gamma + l \delta \leq \frac{r-3}{2} \kappa_r (l/k).\]
It is easy then to check that for some large multiple of the function 
$\alpha_r$ the conditions of Lemma 4.1 are satisfied. 

The first step in formalising this argument is to define an appropriate 
collection of functions $\alpha_r$, which we do as follows:

\demo{Notation}
Let $r \geq 4$ be a positive integer. We write $\beta_r: [0, 
1] \rightarrow [0, \infty)$ for the polynomial function given by
\[\beta_r (z) = 6 z^5 - 15 z^4 + 10 z^3,\]
and $\alpha_r: [o, \infty) \rightarrow [0, \infty)$ for the 
function given by
\[\alpha_r (x) = \left\{\begin{array}{ll}
0 & \mbox{ if $0 \leq x \leq 1 - \frac{1}{2r}$;}\\
\frac{r-4}{4} \beta_r (2rx - (2r-1)) & \mbox{ if $1 - \frac{1}{2r} \leq x 
\leq 1$;}\\
\alpha_r (1/x) & \mbox{ if $x \geq 1$.}\end{array} \right. \]
\Enddemo
 
This slightly bizarre looking set of functions is chosen just so as to 
satisfy the following simple lemma:

\begin{lemma}
For all $r \geq 4$, $\alpha_r$ is a twice-differentiable function such 
that:
\begin{description}
\item{1.}
for $0 \leq x \leq 1$, $0 \leq \alpha_r (x) \leq \frac{r-4}{2} x$;
\item{2.}
$|\alpha'_r (x)| \leq r^2$ and $|\alpha''_r (x)| \leq 20 r^3$ for all $x$.
\end{description}
\end{lemma}

Before we start into the next lemma, we will again need some notation:

\demo{Notation} Suppose that $r \geq 4$ is a fixed positive integer. 
We then write
\[\phi_{r} (k,l) = \alpha_r (l/k) \log (k+l).\]
\Enddemo

Our aim now is to show that $f_r = \exp (-\phi_r)$ (or rather some large 
multiple of it) is an admissible function. The first step towards this is 
contained in the following lemma (note that this is essentially the same 
as Lemma 4 in \cite{T88}):

\begin{lemma} For $k$ and $l$ greater than or equal to $200 r^{10}$, the 
inequalities  
\[\exp\{\phi_r(k,l) - \phi_r(k-m, l)\} \leq 1 + m \Gamma\]
and
\[\exp\{\phi_r(k,l) - \phi_r(k, l-m)\} \leq 1 + m \Delta,\]
where
\[\Gamma = \alpha_r (l/k) \frac{1}{k+l} - \alpha'_r (l/k) \frac{l \log 
(k+l)}{k^2} + \frac{1}{4(k+l)}\]
\[\Delta = \alpha_r (l/k) \frac{1}{k+l} + \alpha'_r (l/k) \frac{\log 
(k+l)}{k} + \frac{1}{4(k+l)}\]
hold for $m = 1, 2$ and $r-1$.
\end{lemma}

\Proof
If we regard $\phi_r(k,l)$ as a function of $k$ with $l$ fixed, then we 
have, using Taylor's Theorem and the fact that $\phi_r$ is twice 
differentiable, that
\[\phi_r(k,l) - \phi_r(k-m,l) = m \frac{\partial\phi_r}{\partial k}(k,l) 
- \frac{m^2}{2} \frac{\partial^2\phi_r}{\partial k^2}(k- \theta m,l)\]
for some $\theta$ between 0 and 1. Now we have that
\[\frac{\partial\phi_r}{\partial k} = \alpha_r (l/k)\frac{1}{k+l} - 
\alpha'_r (l/k) \frac{l \log (k+l)}{k^2},\]
and
\[\frac{\partial^2\phi_r}{\partial k^2}(k,l) = - \alpha_r (l/k) 
\frac{1}{(k+l)^2} - 2 \alpha'_r (l/k) \frac{l}{k^2(k+l)} + 2 \alpha'_r 
(l/k) \frac{l \log (k+l)}{k^3}\]
\[+ \alpha''_r(l/k) \frac{l^2 \log (k+l)}{k^4}.\]
Now note (by using part 2 of Lemma 5.1) that, for $k$ and $l$ both greater 
than or equal to $200 r^{10}$, 
$|\frac{\partial^2\phi_r}{\partial k^2}(k,l)|$ is less than or equal to 
$\frac{1}{4r(k+l)}$. 

Therefore, in this case, we have that
\[\phi_r(k,l) - \phi_r(k-m,l) \leq m \left( \alpha_r (l/k)\frac{1}{k+l} - 
\alpha'_r (l/k) \frac{l \log (k+l)}{k^2} + \frac{1}{8(k+l)}\right).\]
For brevity let's call the right hand side $mx$, noting that for $k$ and 
$l$ greater than or equal to $200 r^{10}$, $mx \leq 1$.

Therefore, using the fact that $e^z \leq 1 + z + z^2$ for $|z| \leq 1$, 
we see that 
\[\exp\{\phi_r(k,l) - \phi_r(k-m,l)\} \leq 1 + mx + m^2 x^2.\]
Note then that, as for the second derivative, by taking $k$ and $l$ 
larger than $200 r^{10}$, we can make $r x^2$ smaller than 
$\frac{1}{8(k+l)}$. Therefore, adding everything together, we see that
\[x + mx^2 \leq \alpha_r (l/k) \frac{1}{k+l} - \alpha'_r (l/k) \frac{l 
\log (k+l)}{k^2} + \frac{1}{4(k+l)},\]
which yields the required result. The result follows similarly for $l$.
\Endproof

We are now ready to tie together everything we have learned in the 
preceding sections to prove a theorem improving the general upper bound 
for Ramsey numbers. This theorem is as follows:

\begin{theorem}
Let $r \geq 4$ be a fixed positive integer. Then there exists a constant 
$c$ such that 
\[r(k+1, l+1) \leq r^{cr^2} \exp\{-\phi_r (k,l)\} \binom{k+l}{k}.\]
\end{theorem}

\Proof
To begin let us suppose that $f_r$ is a function of the form $f_r(a,b) = 
C \exp\{-\phi_r (a,b)\}$ for some fixed constant $C$, and let $n = 
\lfloor f_r(k,l) \binom{k+l}{k} \rfloor = f_r^* (k,l) \binom{k=l}{k}$, 
say. Suppose also that $\kappa_r: [0,\infty) \rightarrow [0, \infty)$ is 
the function given by 
\[\kappa_r (x) = \left\{ \begin{array}{ll}
0 & \mbox{ if $0 \leq x < 1 - \frac{1}{r}$;}\\
\frac{r-3}{2} x & \mbox{ if $1 - \frac{1}{r} \leq x \leq 1$;}\\
\kappa_r(1/x) & \mbox{ if $x \geq 1$.} \end{array} \right.\]
Then, by Lemma 5.1, since $\alpha_r (x) \leq \frac{r-4}{2} x$, we see 
that, for $1 \geq x \geq 1 - \frac{1}{r}$, we have
\[\kappa_r(x) \geq \alpha_r (x) + \frac{1}{2}.\]
If we now choose $k$ and $l$ to both be greater than $200r^{10}$, we can 
apply Lemma 5.2 to see that 
\[\frac{f_r(k-m,l)}{f_r(k,l)} = \exp\{\phi_r(k,l) - \phi_r(k-m, l)\} \leq 
1 + m \Gamma,\]
where
\[\Gamma \leq \alpha_r (l/k) \frac{1}{k+l} - \alpha'_r (l/k) \frac{l \log 
(k+l)}{k^2} + \frac{1}{4(k+l)}.\]
Furthermore, we have that 
\[\frac{f_r(k-m,l)}{f_r^*(k,l)} \leq \left( 1 + \frac{1}{n} \right) 
\frac{f_r(k-m,l)}{f_r(k,l)} \leq 1 + m \gamma,\]
where 
\[\gamma \leq \alpha_r (l/k) \frac{1}{k+l} - \alpha'_r (l/k) \frac{l \log 
(k+l)}{k^2} + \frac{1}{2(k+l)}.\]
Similarly, we have that, for $k$ and $l$ both larger than $200r^{10}$,
\[\frac{f_r(k,l-m)}{f_r^*(k,l)} \leq \left( 1 + \frac{1}{n} \right) 
\frac{f_r(k,l-m)}{f_r(k,l)} \leq 1 + m \delta,\]
where
\[\delta \leq \alpha_r (l/k) \frac{1}{k+l} + \alpha'_r (l/k) \frac{\log 
(k+l)}{k} + \frac{1}{2(k+l)}.\]
Note therefore that 
\[k \gamma + l \delta \leq \alpha_r (l/k) + \frac{1}{2} \leq \kappa_r 
(l/k),\]
provided that $1 \geq \min(\frac{l}{k}, \frac{k}{l}) \geq 1 - 
\frac{1}{r}$. For $\min(\frac{l}{k}, \frac{k}{l}) < 1 - \frac{1}{r}$, we 
have, provided $k$ and $l$ are large (again $200 r^{10}$ will easily 
suffice), that $f_r(a,b)$ is constant (equal to 1) close to $(a,b) = 
(k,l)$ and so we again have
\[k \gamma + l \delta \leq \kappa_r(l/k).\]  
Finally, choose $k$ and $l$ to be sufficiently large, greater than 
$r^{cr}$, for some appropriate $c$, such that Lemma 4.1 holds in the 
following form:
suppose that for $m=1, 2$ and $r-1$, each of the inequalities
\[r(k+1 -m, l+1) \leq f(k-m, l) \binom{k-m+l}{k-m},\]
\[r(k+1, l+1-m) \leq f(k, l-m) \binom{k+l-m}{k},\]
\[\frac{f(k-m,l)}{f^*(k,l)} \leq 1 + m \gamma \mbox{ and }
\frac{f(k,l-m)}{f^*(k,l)} \leq 1 + m \delta\]
holds. Suppose also that $|\gamma|$ and $|\delta|$ are both smaller than 
$r \frac{\log k}{k}$ and that $f(k,l) \geq \exp \{-r \frac{l}{k} \log k \}$.
Then, provided that 
\[k \gamma + l \delta \leq \kappa_r(l/k),\]
we have that
\[r(k+1,l+1) \leq f(k,l) \binom{k+l}{k}.\]

To conclude, suppose that $N > \max(200 r^{10}, r^{cr}) 
= r^{cr}$, for $c$ chosen large enough, and consider the function 
$f_r(a,b) = (2N)^r \exp\{ -\phi_r(a,b)\}$. For either $a$ or $b$ less 
than or equal to $N$ we have straightforwardly that for $a \geq b$ with 
$b \leq N$,   
\[f_r(a,b) \geq \frac{(2N)^r}{(a+b)^{rb/a}} \geq 1,\]  
using the fact that $(a+b)^{b/a}$ is a decreasing function in $a$.
Now, both $\gamma$ and $\delta$ defined above are less than or equal to 
$r \frac{\log k}{k}$, and $f(k,l)$ is certainly larger than $\exp \{-r 
\frac{l}{k} \log k \}$. Finally, we have by the construction of $\phi_r$ 
and the choice of $N$ that 
\[k \gamma + l \delta \leq \kappa_r(l/k).\]
Consequently our induction holds good with this function $f_r$.
\Endproof

Theorem 1.2 is now a straightforward consequence of this theorem. The 
simple proof is in fact contained in the following proof of Theorem 1.1:

\demo{Proof of Theorem 1.1}
From Theorem 5.1, we know that, for integers $r \geq 5$, 
\begin{eqnarray*}
r(k+1, k+1) & \leq & r^{cr^2} \exp\{-\phi_r (k,k)\} \binom{2k}{k}\\
& \leq & r^{cr^2} \exp\{ -\frac{r-4}{4} \log (2k)\} \binom{2k}{k} \\
& \leq & \frac{r^{cr^2}}{k^{dr}} \binom{2k}{k}
\end{eqnarray*}
for some fixed constants $c$ and $d$.

If now, for any sufficiently large $k$, we take $r = \lfloor \frac{d \log 
k}{2c \log \log k} \rfloor$ (a value which is close to that which 
minimises $\frac{r^{cr^2}}{k^{dr}}$), we see that for some constant $C$ we have 
\[r(k+1, k+1) \leq k^{- C \frac{\log k}{\log \log k}} \binom{2k}{k},\]
as required.
\Endproof

\demo{Acknowledgements}
I would like to thank B\'{e}la Bollob\'{a}s, Tim Gowers, Ben Green and Tom 
Sanders for their advice and encouragement. 
\Enddemo


\references{}

\bibitem[CGW89]{CGW89}
{F.R.K. Chung, R.L. Graham, R.M. Wilson,}
{Quasi-random graphs,}
{\it Combinatorica} {\bf 9} {(1989), no. 4, 345-362.}

\bibitem[ES35]{ES35}
{P. Erd\H{o}s, G. Szekeres,}
{A combinatorial problem in geometry,}
{\it Compositio Math.} {\bf 2} {(1935), 463-470.}

\bibitem[GR87]{GR87}
{R. L. Graham, V. R\"{o}dl,}
{Numbers in ramsey theory,}
{\it Surveys in Combinatorics, London Math. Soc. Lecture Note Series no. 
123, Cambridge Univ. Press, London} {(1987), 111-153}

\bibitem[R30]{R30}
{F. P. Ramsey,}
{On a problem of formal logic,}
{\it Proc. London Math. Soc.} {\bf 30} {(1930), 264-286}

\bibitem[T87]{T87}
{A. Thomason,}
{Pseudorandom graphs,}
{\it Random Graphs '85 (Pozna\'{n}, 1985), North-Holland Math. Stud., 
vol. 144, North-Holland, Amsterdam-New York,} {(1987), 307-331.}

\bibitem[T88]{T88}
{A. Thomason,}
{An upper bound for some ramsey numbers,}
{\it J. Graph Theory} {\bf 12} {(1988), 509-517}

\Endrefs
\end{document}